\documentclass[reqno,12]{amsart}
\usepackage{amssymb}
\usepackage{amsmath, amsthm}
\usepackage{hyperref}
\usepackage{enumerate}
\usepackage{graphicx}
\usepackage{bm}

\AtBeginDocument{\noindent\small
\vspace{9mm}}
\begin{document}
\title[On a new identity for the H-function with applications $\cdots$]{On a new identity for the H-function with applications to the summation of hypergeometric series}
\author[Arjun K. Rathie, L.C.S.M. Ozelim  and P.N. Rathie ]{Arjun K. Rathie, L.C.S.M. Ozelim  and P.N. Rathie }
\address{$^1$ Department of Mathematics,  Central University of Kerala,  Kasaragod- 671316, Kerala, India}
\email{akrathie@cukerala.ac.in}
\address{$^2$ Dept. of Civil and Environmental Engineering, University of Brasilia, Brasilia- DF, 70910-900, Brazil}
\email{luanoz@gmail.com} 
\address{$^3$ Dept. of Statistics and Applied Mathematics, Federal University of Ceara, Fortaleza- CE, 60440-900, Brazil}
\begin{abstract}
Using generalized hypergeometric functions to perform symbolic manipulation of equations is of great importance to pure and applied scientists. There are in the literature a great number of identities for the Meijer-G function. On the other hand, when more complex expressions arise, the latter function is not capable of representing them. The H-function is an alternative to overcome this issue, as it is a generalization of the Meijer-G function. In the present paper, a new identity for the H-function is derived. In short, this result enables one to split a particular H-function into the sum of two other H-functions. The new relation in addition to an old result are applied to the summation of hypergeometric series. Finally, some relations between H-functions and elementary functions are built. 
\end{abstract}
\maketitle
\newtheorem{theorem}{Theorem}[section]
\newtheorem{lemma}[theorem]{Lemma}
\newtheorem{proposition}[theorem]{Proposition}
\newtheorem{corollary}[theorem]{Corollary}
\newtheorem{Definition}[theorem]{Definition}
\newtheorem*{remark}{Remark}
\newtheorem{integral}[theorem]{Integral}
\section{Introduction}
Special functions have proven to be essential tools while dealing with the formal mathematical manipulation of equations. In fact, most of the computational softwares which perform symbolic operations consider generalized hypergeometric functions to do so.

Generalized hypergeometric functions of the type $_{p}F_{q}$ have been extensively studied. For example, in the works of \cite{Bailey,HMS2012,Choi,HMS2014}, a series of identities have been derived for this function.

These functions, on the other hand, are able to represent just a small share of the mathematical relations commonly considered in science. Thus, more general hypergeometric functions must be considered. This is the case of the Meijer-G function \cite{Mathai2010}.

Mathematica software, for example, vastly relies on the Meijer-G function to perform integration, differentiation and algebric manipulation of standard and special functions. This comes from the fact that most of the functions which are used in science are representable in terms of this special function.

When more complex expressions arise, Meijer-G functions are not capable of representing the functional relations which show up. Thus, a more general function is needed for this task. This is where the H-function, which is a generalization of the Meijer-G function, can be used. The H-function is a powerful hypergeometric function whose importance in pure and applied sciences has been considerably discussed \cite{Mathai2010, Springer1979}.

Useful identities for this special function have been presented in \cite{Mathai2010,AKR1981}. In the present paper, a new identity for the H-function has been derived. It is shown how this identitiy can be used to provide closed form representations for hypergeometric summations. Besdies, new relations between H-functions and elementary functions are presented.

In order to better familiarize the reader with this special function, the next section presents some basic concepts regarding the H-function.

\section{H-function}

The H - function (see \cite{Mathai2010} ) is defined, as a contour
complex integral which contains gamma functions in their
integrands, by
\begin{gather}
H_{p,q}^{m,n}\left[z\;\bigg|\begin{array}{cccccc}
(a_{1},A), & \ldots, & (a_{n},A_{n}), & (a_{n+1},A_{n+1}), & \ldots, & (a_{p},A_{p})\\
(b_{1},B_{1}), & \ldots, & (b_{m},B_{m}), & (b_{m+1},B_{m+1}), &
\ldots, & (b_{q},B_{q})
\end{array}\right]\nonumber \\
=\frac{1}{2\pi i}\int_{L}\frac{{\displaystyle
\prod_{j=1}^{m}\Gamma(b_{j}+B_{j}s){\displaystyle
\prod_{j=1}^{n}\Gamma(1-a_{j}-A_{j}s)}}}{{\displaystyle
\prod_{j=m+1}^{q}\Gamma(1-b_{j}-B_{j}s){\displaystyle
\prod_{j=n+1}^{p}\Gamma(a_{j}+A_{j}s)}}}z^{-s}ds,\label{eq:Hfunccontour}
\end{gather}
where \ $A_{j}$ and $B_{j}$ are assumed to be positive quantities
and all the $a_{j}$ and $b_{j}$ may be complex. The contour $L$ runs
from $c-i\infty$ to $c+i\infty$ such that the poles of
$\Gamma(b_{j}+B_{j}s)$, \ $j=1,\ldots,m$ lie to the left of $L$
and the poles of $\Gamma(1-a_{j}-A_{j}s)$, $j=1,\ldots,n$ \ lie to
the right of $L$.

By performing the variable change $s \to -r$ and adjusting the contour $L$ to $L^{*}$, where the integral runs from $c^{*} - i \infty$ to  $c^{*} + i \infty$, the H-function can be alternatively defined as:

\begin{gather}\label{eq:Hfunccontourneg}
H_{p,q}^{m,n}\left[z\;\bigg|\begin{array}{cccccc}
(a_{1},A), & \ldots, & (a_{n},A_{n}), & (a_{n+1},A_{n+1}), & \ldots, & (a_{p},A_{p})\\
(b_{1},B_{1}), & \ldots, & (b_{m},B_{m}), & (b_{m+1},B_{m+1}), &
\ldots, & (b_{q},B_{q})
\end{array}\right]\nonumber \\
=\frac{1}{2\pi i}\int_{L^{*}}\frac{{\displaystyle
\prod_{j=1}^{m}\Gamma(b_{j}-B_{j}r){\displaystyle
\prod_{j=1}^{n}\Gamma(1-a_{j}+A_{j}r)}}}{{\displaystyle
\prod_{j=m+1}^{q}\Gamma(1-b_{j}+B_{j}r){\displaystyle
\prod_{j=n+1}^{p}\Gamma(a_{j}-A_{j}r)}}}z^{r}dr,
\end{gather}

\noindent for which the same parameter domains' restrictions apply.

By considering the definition in (\ref{eq:Hfunccontourneg}), the H-function can be expressed in computable form as \cite{Mathai2010}:

When the poles of ${\displaystyle \prod_{j=1}^{m}\Gamma(b_{j}-B_{j}r)}$
are simple, we have:

\begin{align} \label{hcomput1}
H_{p\;\; q}^{m\; n}(z) & =\sum_{h=1}^{m}\sum_{\nu=0}^{\infty}\frac{{\displaystyle \prod_{j=1\neq h}^{m}\Gamma\left(b_{j}-B_{j}\frac{b_{h}+\nu}{B_{h}}\right)}}{{\displaystyle \prod_{j=m+1}^{q}\Gamma\left(1-b_{j}+B_{j}\frac{b_{h}+\nu}{B_{h}}\right)}}\times\nonumber \\
 & \quad\times\frac{{\displaystyle \prod_{j=1}^{n}\Gamma\left(1-a_{j}+A_{j}\frac{b_{h}+\nu}{B_{j}}\right)}}{{\displaystyle \prod_{j=n+1}^{p}\Gamma\left(a_{j}-A_{j}\frac{b_{h}+\nu}{B_{h}}\right)}}\frac{(-1)^{\nu}z^{(b_{h}+\nu)/B_{h}}}{\nu!B_{h}}
\end{align}

\noindent for $z\neq0$ if \ $\delta>0$ \ and \ for \ $0<|z|<D^{-1}$ if
\ $\delta=0$, where \ $\delta=\sum_{j=1}^{p}B_{j}-\sum_{j=1}^{q}A_{j}$
and \ $D=\prod_{j=1}^{p}A_{j}^{A_{j}}/\prod_{j=1}^{q}B_{j}^{B_{j}}$.

When the poles of ${\displaystyle \prod_{j=1}^{n}\Gamma(1-a_{j}+A_{j}r)}$
are simple, we have 

\begin{align} \label{hcomput2}
H_{p\; q}^{m\; n}(z) & =\sum_{h=1}^{n}\sum_{\nu=0}^{\infty}\frac{{\displaystyle \prod_{j=1\neq h}^{n}\Gamma\left(1-a_{j}-A_{j}\frac{1-a_{h}+\nu}{A_{h}}\right)}}{{\displaystyle \prod_{j=n+1}^{p}\Gamma\left(a_{j}+A_{j}\frac{1-a_{h}+\nu}{A_{h}}\right)}}\times\nonumber \\
 & \quad\times\frac{{\displaystyle \prod_{j=1}^{m}\Gamma\left(b_{j}+B_{j}\frac{1-a_{h}+\nu}{A_{h}}\right)}}{{\displaystyle \prod_{j=m+1}^{q}\Gamma\left(1-b_{j}-B_{j}\frac{1-a_{h}+\nu}{A_{h}}\right)}}\frac{(-1)^{\nu}(1/z)^{(1-a_{h}+\nu)/A_{h}}}{\nu!A_{h}}
\end{align}

\noindent for $z\neq0$ if $\delta<0$ and for $|z|>D^{-1}$ if $\delta=0$.

Both representations above apply when the poles of the gamma function
in the numerator of the quotients are simple. When this simplification
does not hold, residue theorem has to be applied. For details about
this theorem, one may refer to \cite{Springer1979}.

Another hypergeometric function which is of interest in the present paper is the ${}_{p}F_{q}$, defined as :
\begin{equation}\label{pfq}
{}_{p}F_{q} \left[\begin{array}{c} a_1, \ldots, a_p \\b_1, \ldots, b_q \end{array} \; ; z \right] = \sum_{n=0}^{\infty} \frac{(a_1)_n \ldots (a_p)_n}{(b_1)_n \ldots (b_q)_n} \frac{x^n}{n!}
\end{equation}
\noindent where the symbols follow the same constraints as in the case of the H-function. Also, $(a)_n$ denotes the Pochhammer symbol, which can be defined in termos of the Gamma function as:
\begin{equation}\label{poch}
(a)_n = \frac{\Gamma (a+n)}{\Gamma(a)}
\end{equation}
Both the H-function and the ${}_{p} F_{q}$ function may be related by the following formula:
\begin{equation}\label{pfqh}
{}_{p} F_{q} \left[ \begin{array}{c}a_1, \ldots, a_p \\b_1, \ldots, b_q \end{array} \; ; z \right] = \frac{\prod_{k=1}^{q} \Gamma (b_k)}{\prod_{k=1}^{p} \Gamma (a_k)} H_{p,q+1}^{1,p}\left[-z|\begin{array}{c} (1-a_1,1),\ldots, (1-a_p,1)\\ (0,1),(1-b_1,1),\ldots, (1-b_q,1)\end{array}\right] 
\end{equation}
\section{Identities presented in the literature}
In the present section, a few identities presented in the literature are shown in order to better familiarize the reader with the mathematics behind the proofs \cite{Mathai2010,Brychkov,Prudnikov}.
\begin{equation} \label{mult}
\prod_{j=0}^{k-1} \Gamma \left( z + \frac{j}{k} \right) = \Gamma (k z) (2 \pi)^{\frac{k-1}{2}} k^{\frac{1}{2} - k z}
\end{equation}
\noindent where $k$ is a positive integer; $kz \in \mathbb{C} \backslash Z_{0}^{-}$.
\begin{equation} \label{cosg}
cos (\pi z) = \frac{\pi}{\Gamma \left( \frac{1}{2} + z \right) \Gamma \left( \frac{1}{2} - z \right) }
\end{equation} 
\begin{equation} \label{cose}
cos (\pi z) = \frac{e^{i \pi z} + e^{- i \pi z}}{2}
\end{equation} 
\begin{equation}
H_{p,q}^{m,n}\left[z|\begin{array}{c}
(a_p,A_p)\\
(b_q,B_q)
\end{array}\right] =
kH_{p,q}^{m,n}\left[z^k|\begin{array}{c}
(a_p,kA_p)\\
(b_q,kB_q)
\end{array}\right], \quad k>0.
\label{hexp}
\end{equation}

\begin{equation} \label{exph}
H_{0,1}^{1,0}\left[z\;\bigg|\begin{array}{c}
- \\
(0,1)
\end{array}\right] = e^{-z}
\end{equation}

\begin{equation} \label{h2}
H_{1,2}^{1,1}\left[-z\;\bigg|\begin{array}{c}
(0,1)\\
(0,1), (-1,1)
\end{array}\right] = \frac{e^z-1}{z}
\end{equation}

\begin{equation} \label{h3}
H_{2,3}^{1,2}\left[-z\;\bigg|\begin{array}{c}
(0,1),(-1,1)\\
(0,1),(-1,1),(-2,1)
\end{array}\right] = \frac{e^z-1-z}{z^2}
\end{equation}

In \cite{AKR1981}, an interesting relation was derived to split an H-function into the sum of two other H-functions. This relation can be expressed as \cite{AKR1981}:
Let $z \in \mathbb{C}$, then:
\begin{align} \label{eq:HfuncIDAKR}
H_{p+1,q+1}^{m,n}& \left[z\;\bigg|\begin{array}{c}
(a_{1},A_{1}), \ldots, (a_{p},A_{p}),(\alpha,\lambda)\\
(b_{1},B_{1}), \ldots, (b_{q},B_{q}),(\alpha,\lambda)   \end{array}\right]  \nonumber \\
&  = \frac{1}{2 \pi i} \left( e^{i \pi \alpha} H_{p,q}^{m,n}\left[e^{- i \pi \lambda} z\;\bigg|\begin{array}{c}(a_{1},A_{1}), \ldots, (a_{p},A_{p})\\ (b_{1},B_{1}), \ldots, (b_{q},B_{q})
\end{array}\right]  \right.  \nonumber \\
& \left. \qquad \qquad   -  e^{- i \pi \alpha} H_{p,q}^{m,n}\left[e^{i \pi \lambda} z\;\bigg|\begin{array}{c}
(a_{1},A_{1}), \ldots, (a_{p},A_{p})\\ (b_{1},B_{1}), \ldots, (b_{q},B_{q}) \end{array}\right] \right), 
\end{align}
In the present paper, an alternative splitting relation is derived, as shall be seen in the next section.

\section{Results}

At first, one identity for the H-function is presented. Then, applications of the new identity derived are shown together with (\ref{eq:HfuncIDAKR}).

\begin{theorem}

Let $z \in \mathbb{C}$, then:
\begin{align}
H_{p,q}^{m,n} & \left[z\;\bigg|\begin{array}{c}
(\alpha,\lambda), (a_{2},A_{2}), \ldots, (a_{p},A_{p}) \\ (\alpha,\lambda), (b_{2},B_{2}), \ldots, (b_{q},B_{q}) \end{array}\right] \nonumber \\
& = e^{i \pi \alpha} H_{p,q}^{m,n}\left[e^{- i \pi \lambda} z\;\bigg|\begin{array}{c}
(2\alpha,2\lambda), (a_{2},A_{2}), \ldots, (a_{p},A_{p})\\ (2\alpha,2\lambda), (b_{2},B_{2}), \ldots, (b_{q},B_{q})\end{array}\right] \nonumber \\
& + e^{- i \pi \alpha} H_{p,q}^{m,n}\left[e^{i \pi \lambda} z\;\bigg|\begin{array}{c}
(2\alpha,2\lambda), (a_{2},A_{2}), \ldots, (a_{p},A_{p})\\
(2\alpha,2\lambda), (b_{2},B_{2}), \ldots, (b_{q},B_{q})
\end{array}\right],\label{eq:HfuncID}
\end{align}

\begin{proof}
At first, by using (\ref{eq:Hfunccontour}), one shall consider the contour integral representation of the H-function in (\ref{eq:HfuncID}), given as:

\begin{align} \label{split1}
H_{p,q}^{m,n} & \left[z\;\bigg|\begin{array}{c}(\alpha,\lambda), \ldots, (a_{p},A_{p})\\(\alpha,\lambda), \ldots, (b_{q},B_{q}) \end{array}\right]  \nonumber \\
 & = \frac{1}{2\pi i}\int_{L}\frac{\Gamma(\alpha + \lambda s) \Gamma(1 - \alpha - \lambda s) {\displaystyle \prod_{j=2}^{m}\Gamma(b_{j}+B_{j}s){\displaystyle \prod_{j=2}^{n}\Gamma(1-a_{j}-A_{j}s)}}}{{\displaystyle \prod_{j=m+1}^{q}\Gamma(1-b_{j}-B_{j}s){\displaystyle \prod_{j=n+1}^{p}\Gamma(a_{j}+A_{j}s)}}} z^{-s}ds
\end{align}

By considering (\ref{mult}), it is easy to see that:

\begin{eqnarray} \label{split2}
\Gamma(\alpha + \lambda s) \Gamma \left(\frac{1}{2} +\alpha + \lambda s \right) = \Gamma (2 \alpha + 2 \lambda s) (2 \pi)^{\frac{1}{2}} 2^{\frac{1}{2} - 2 \alpha - 2 \lambda s} \nonumber \\
\Gamma \left(\frac{1}{2} -\alpha - \lambda s \right) \Gamma \left(1 -\alpha - \lambda s \right) = \Gamma (1- 2 \alpha - 2 \lambda s) (2 \pi)^{\frac{1}{2}} 2^{\frac{1}{2} -1+ 2 \alpha + 2 \lambda s} 
\end{eqnarray}

Also, (\ref{split2}) implies that:

\begin{equation} \label{split3}
\Gamma(\alpha + \lambda s) \Gamma \left(1 -\alpha - \lambda s \right) =2 \pi \frac{\Gamma (2 \alpha + 2 \lambda s) \Gamma (1- 2 \alpha - 2 \lambda s) }{\Gamma \left(\frac{1}{2} +\alpha + \lambda s \right) \Gamma \left(\frac{1}{2} -\alpha - \lambda s \right)}
\end{equation} 

Equation (\ref{split3}) can be further simplified by using (\ref{cosg}) and (\ref{cose}), resulting in:

\begin{equation} \label{split4}
\Gamma(\alpha + \lambda s) \Gamma \left(1 -\alpha - \lambda s \right) = \Gamma (2 \alpha + 2 \lambda s) \Gamma (1- 2 \alpha - 2 \lambda s) (e^{i \pi (\alpha + \lambda s)} + e^{- i \pi (\alpha + \lambda s)})
\end{equation} 

Finally, by inserting (\ref{split4}) into (\ref{split1}), (\ref{eq:HfuncID}) is retrieved. 

\end{proof}
\end{theorem}

The following two forms of (\ref{eq:HfuncID}) and (\ref{eq:HfuncIDAKR}) are easily derived:

\begin{itemize}
\item By taking $(a_1,A_1)=(b_1,B_1)=(\alpha,\lambda)$ in (\ref{eq:HfuncIDAKR}), one gets:
\begin{align} \label{eq:HfuncIDAKR2}
H_{p-1,q-1}^{m-1,n-1} & \left[z\;\bigg|\begin{array}{c}(a_{2},A_{2}), \ldots, (a_{p},A_{p})\\
(b_{2},B_{2}), \ldots, (b_{q},B_{q})\end{array}\right] \nonumber\\
& = \frac{1}{2 \pi i} \bigg( e^{i \pi \alpha} H_{p,q}^{m,n}\left[e^{- i \pi \lambda} z\;\bigg|\begin{array}{c}(\alpha,\lambda),(a_{2},A_{2}), \ldots, (a_{p},A_{p})\\ (\alpha,\lambda),(b_{2},B_{2}), \ldots, (b_{q},B_{q})
\end{array}\right] \nonumber \\
& \quad - e^{- i \pi \alpha} H_{p,q}^{m,n}\left[e^{i \pi \lambda} z\;\bigg|\begin{array}{c}
(\alpha,\lambda),(a_{2},A_{2}), \ldots, (a_{p},A_{p})\\ (\alpha,\lambda),(b_{2},B_{2}), \ldots, (b_{q},B_{q})\end{array}\right] \bigg),
\end{align}
\noindent for $p \geq n \geq 1$ and $q \geq m \geq 1$.
\item Similarly, taking $(a_p,A_p) = (b_q,B_q) = (\alpha,\lambda)$ in (\ref{eq:HfuncID}), one obtains:
\begin{align}
& H_{p-2,q-2}^{m-1,n-1}  \left[z\;\bigg|\begin{array}{c}(a_{2},A_{2}), \ldots, (a_{p-1},A_{p-1})\\
(b_{2},B_{2}), \ldots, (b_{q-1},B_{q-1})\end{array}\right] \nonumber\\
& = e^{i \pi \alpha} H_{p,q}^{m,n}\left[e^{- i \pi \lambda} z\;\bigg|\begin{array}{c}
(2\alpha,2\lambda), (a_{2},A_{2}), \ldots, (a_{p-1},A_{p-1}),(\alpha,\lambda)\\
(2\alpha,2\lambda), (b_{2},B_{2}), \ldots, (b_{q-1},B_{q-1}),(\alpha,\lambda)
\end{array}\right] \nonumber \\
& \quad + e^{- i \pi \alpha} H_{p,q}^{m,n}\left[e^{i \pi \lambda} z\;\bigg|\begin{array}{c}
(2\alpha,2\lambda), (a_{2},A_{2}), \ldots, (a_{p-1},A_{p-1}),(\alpha,\lambda)\\
(2\alpha,2\lambda), (b_{2},B_{2}), \ldots, (b_{q-1},B_{q-1}),(\alpha,\lambda)
\end{array}\right],\label{eq:HfuncIDn2}
\end{align}

\noindent for $p-1 \geq n \geq 1$ and $q-1 \geq m \geq 1$.

\end{itemize}

Equations (\ref{eq:HfuncIDAKR}) and (\ref{eq:HfuncID}) may be applied to split H-functions into the sum of H-function. In the present paper, an interesting application to the summation of hypergeometric series is described. At first, a general summation may be studied. Then, special cases are discussed case by case.

\begin{theorem}

Let $x \in \mathbb{R}$, then:

\begin{equation} \label{refhyp}
\sum_{n=0}^{\infty} \frac{x^{\alpha n} \beta^{n}}{\Gamma (\gamma n + \delta +1)} = \frac{1}{\Gamma (\delta +1)} {}_{1} F_{\gamma} \left[ \begin{array}{c}
1 \\
\frac{\delta + 1}{\gamma},\frac{\delta + 2}{\gamma}, \ldots, \frac{\delta + \gamma}{\gamma}
\end{array} \; ; \beta \frac{x^{\alpha}}{\gamma^{\gamma}} \right]
\end{equation}

The left hand side of (\ref{refhyp}) is a generalization of the Mittag-Leffler function \cite{Mathai2010,Brychkov}.

\begin{proof}

The summation in (\ref{refhyp}) can be expressed in terms of the hypergeometric function ${}_{p} F_{q}$ by noticing that:

\begin{equation} \label{expan}
\Gamma (\gamma n + \delta + 1) = \frac{\prod_{k=0}^{\gamma -1} \Gamma \left( n + \frac{\delta + 1 + k}{\gamma} \right)}{(2 \pi)^{\frac{\gamma -1}{2}} \gamma^{-\frac{1}{2} - \gamma n - \delta}}
\end{equation}

By combining equations (\ref{refhyp}) and (\ref{expan}), the former may be rewritten as:

\begin{equation} \label{hyp2}
\frac{1}{\Gamma (\delta +1)}  \sum_{n=0}^{\infty} \frac{(\beta x^{\alpha} \gamma^{-\gamma})^{n} (1)_n}{\prod_{k=0}^{\gamma-1} \left(\frac{\delta + 1 +k}{\gamma}\right)_n n!}
\end{equation}

By the representation in (\ref{pfq}), equation (\ref{hyp2}) implies (\ref{refhyp}).

\end{proof}
\end{theorem}

\begin{corollary}

The summation in (\ref{refhyp}) can be expressed in terms of the H-function as:

\begin{equation} \label{refhyp2}
\sum_{n=0}^{\infty} \frac{x^{\alpha n} \beta^{n}}{\Gamma(\gamma n + \delta +1)} =  
H_{1,2}^{1,1}\left[-\beta x^{\alpha}|\begin{array}{c}
(0,1)\\
(0,1),(-\delta,\gamma)
\end{array}\right] 
\end{equation}

\begin{proof}

By means of (\ref{pfqh}), (\ref{refhyp}) is expressible in terms of the H-function as:

\begin{align} \label{refhyp3}
\sum_{n=0}^{\infty} \frac{x^{\alpha n} \beta^{n}}{\Gamma (\gamma n + \delta +1)} &=  \frac{1}{\Gamma (\delta +1)} \prod_{k=0}^{\gamma -1} \Gamma (\frac{\gamma -\delta - 1- k}{\gamma}) \nonumber \\
& \quad . H_{1,\gamma+1}^{1,1}\left[-\beta \frac{x^{\alpha}}{\gamma^{\gamma}}|\begin{array}{c}
(0,1)\\ (0,1),(1-\frac{\delta + 1}{\gamma},1),\ldots, (1-\frac{\delta + \gamma}{\gamma},1)
\end{array}\right] 
\end{align}

On the other hand, by using the contour integral representation (\ref{eq:Hfunccontour}) of (\ref{refhyp3})  and the multiplication theorem for the Gamma function (\ref{mult}), (\ref{refhyp2}) follows from (\ref{refhyp3}). 

\end{proof}

\end{corollary}

Based on the general results presented, one shall proceed to the applications sections.

\section{Applications}

In the present section, results obtained for the summation of hypergeometric series are presented. These results are important as present interesting relations between generalized hypergeometric functions and elementary functions. 

\subsection{When $\gamma=2$ and $\beta = -1$ in (\ref{refhyp})}

For the cases where $\beta=-1$, the results may be obtained by performing the transformation $x^{\alpha/2} \to i x^{\alpha/2}$. This provides the results in terms of standard trigonometrical functions. Also, these are the cases treated in \cite{HMS2014} for $\Omega_i = 1$, $\forall i \in \mathbb{N}$, where the parameters considered are: $\gamma=2$ ($\delta = 0$ and $1$) and $\gamma=3$ ($\delta = 0, 1$ and $2$).

Another interesting result may be discussed for the case where $\gamma=2$ and $\beta = -1$.

\subsubsection{Case when $\delta = 0$}

At first, let one consider (\ref{refhyp}) for this case:

\begin{equation} \label{hypn}
\sum_{n=0}^{\infty} \frac{x^{\alpha n} (-1)^{n}}{\Gamma (2 n +1)} =
{}_{0} F_{1} \left[ \begin{array}{c}
- \\
\frac{1}{2}
\end{array} \; ; - \frac{x^{\alpha}}{4} \right]
\end{equation}

The literature shows that when $\alpha = 2$, we have \cite[eq.(3.2)]{Bailey}:

\begin{equation} \label{hypnbai1}
e^{x} {}_{0} F_{1} \left[ \begin{array}{c}
- \\
\frac{1}{2}
\end{array} \; ; - \frac{x^{2}}{4} \right] = 
\sum_{n=0}^{\infty} 2^{n/2} cos \left( \frac{n \pi}{4}\right) \frac{x^n}{n!}
\end{equation}

By considering (\ref{cosg}), the right hand side of (\ref{hypnbai1}) may be rewritten as:

\begin{equation} \label{hypnbai2}
e^{x} {}_{0} F_{1} \left[ \begin{array}{c}
- \\
\frac{1}{2}
\end{array}\; ; - \frac{x^{2}}{4} \right] = \pi  
\sum_{n=0}^{\infty} \frac{2^{n/2}}{\Gamma \left( \frac{1}{2} + \frac{n}{4} \right) \Gamma \left( \frac{1}{2} - \frac{n}{4} \right)} \frac{x^n}{n!}
\end{equation}

Finally, by considering the H-function series representation given in (\ref{hcomput1}), (\ref{hypnbai2}) provides:

\begin{equation} \label{hypnbai3}
e^{x} {}_{0} F_{1} \left[\begin{array}{c}
- \\
\frac{1}{2}
\end{array} \; ; - \frac{x^{2}}{4} \right] = \pi  
H_{1,2}^{1,0}\left[-\sqrt{2} x \bigg| \begin{array}{c}
(\frac{1}{2},\frac{1}{4})\\
(0,1),(\frac{1}{2},\frac{1}{4})
\end{array}\right]
\end{equation}

The alternative representation (\ref{refhyp2}) of the left hand side of (\ref{hypnbai3}) results in the following relations:

\begin{equation} \label{hypnbai4}
e^{x} {}_{0} F_{1} \left[\begin{array}{c}
- \\
\frac{1}{2}
\end{array} \; ; - \frac{x^{2}}{4} \right]= e^{x} 
H_{1,2}^{1,1}\left[ x^{2} \bigg| \begin{array}{c}
(0,1)\\
(0,1),(0,2)
\end{array}\right] = \pi  
H_{1,2}^{1,0}\left[-\sqrt{2} x \bigg| \begin{array}{c}
(\frac{1}{2},\frac{1}{4})\\
(0,1),(\frac{1}{2},\frac{1}{4})
\end{array}\right]
\end{equation}

By considering the transformation $x \to i x$, (\ref{hypnbai4}) is also expressible in terms of elementary functions, as $e^x cos (x)$, on using \cite{Brychkov}.

\subsubsection{Case when $\delta = 1$}

For the case where $\alpha = 2$ and $\delta = 1$, by considering the relation in \cite[eq.(3.3)]{Bailey} and following a similar procedure as in the case where $\delta = 0$, the following relation is obtained:

\begin{align} \label{hypnbai5}
e^{x} {}_{0} F_{1} \left[ \begin{array}{c} - \\ \frac{3}{2} \end{array} \; ; - \frac{x^{2}}{4} \right] & = e^{x} H_{1,2}^{1,1}\left[ x^{2} \bigg| \begin{array}{c} (0,1)\\ (0,1),(-1,2) \end{array}\right]\nonumber \\ 
 & = \sqrt{2} \pi H_{2,3}^{1,1}\left[-\sqrt{2} x \bigg| \begin{array}{c} (0,1),\left(\frac{3}{4},\frac{1}{4}\right)\\ (0,1),(-1,1),\left(\frac{3}{4},\frac{1}{4} \right) \end{array}\right]
\end{align}
On the other hand, by considering the transformation $x \to i x$, (\ref{hypnbai5}) can be expresses in terms of trigonometrical functions, as $x^{-1} e^x sin(x)$ on using \cite{Brychkov}.

\subsection{When $\gamma=3$ and $\beta = -1$ in (\ref{refhyp})}

For the next subcases, a similar procedure as in the cases where $\gamma=2$ and $\beta = -1$ can be followed. Thus, by considering \cite[eq.(3.4)]{Bailey}, \cite[eq.(3.5)]{Bailey} and \cite[eq.(3.6)]{Bailey} for the cases where $\delta = 0$, $\delta = 1$ and $\delta = 2$, respectively, the following relations between H-functions are retrieved:

\begin{align} \label{hypnbai6}
e^{x} {}_{0} F_{2} \left[ \begin{array}{c}- \\\frac{1}{3},\frac{2}{3}\end{array} \; ; - \frac{x^{3}}{27} \right] & = e^{x} H_{1,2}^{1,1}\left[ x^{3} \bigg| \begin{array}{c} (0,1)\\
(0,1),(0,3) \end{array}\right]\nonumber \\ 
&  = \frac{1}{3}+ \frac{2 \pi}{3} 
H_{1,2}^{1,0}\left[-\sqrt{3} x \bigg| \begin{array}{c}
(\frac{1}{2},\frac{1}{6})\\
(0,1),(\frac{1}{2},\frac{1}{6})
\end{array}\right]
\end{align}

\begin{align} \label{hypnbai7}
e^{x} {}_{0} F_{2} \left[ \begin{array}{c}
- \\
\frac{2}{3},\frac{4}{3}
\end{array} \; ; - \frac{x^{3}}{27} \right]& = e^{x} 
H_{1,2}^{1,1}\left[ x^{3} \bigg| \begin{array}{c}
(0,1)\\
(0,1),(-1,3)
\end{array}\right] \nonumber \\ 
&  = \frac{2 \pi}{\sqrt{3}}  
H_{2,3}^{1,1}\left[-\sqrt{3} x \bigg| \begin{array}{c}
(0,1),(\frac{2}{3},\frac{1}{6})\\
(0,1),(-1,1),(\frac{2}{3},\frac{1}{6})
\end{array}\right]
\end{align}

\begin{align} \label{hypnbai8}
e^{x} {}_{0} F_{2} \left[ \begin{array}{c}
- \\
\frac{4}{3},\frac{5}{3}
\end{array} \; ; - \frac{x^{3}}{27} \right] & = e^{x} 
H_{1,2}^{1,1}\left[ x^{3} \bigg| \begin{array}{c}
(0,1)\\
(0,1),(-2,3)
\end{array}\right] \nonumber \\ 
&  = 4 \pi  
H_{3,4}^{1,2}\left[-\sqrt{3} x \bigg| \begin{array}{c}
(0,1),(-1,1),(\frac{5}{6},\frac{1}{6})\\
(0,1),(-1,1),(-2,1),(\frac{5}{6},\frac{1}{6})
\end{array}\right]
\end{align}

Each of the equations from (\ref{hypnbai6}) to (\ref{hypnbai8}) can be expressed in terms of elementary functions. This can be accomplished by using (\ref{eq:HfuncIDAKR}) and \cite{Brychkov}.

\subsubsection{Case when $\delta = 0$}

By combining (\ref{hypnbai6}) and (\ref{eq:HfuncIDAKR}), the following is obtained:

\begin{align} \label{hypnbai6-1}
1+ \frac{2 \pi}{3}  
H_{1,2}^{1,0} & \left[-\sqrt{3} x \bigg| \begin{array}{c}
(\frac{1}{2},\frac{1}{6})\\ (0,1),(\frac{1}{2},\frac{1}{6})
\end{array}\right] \nonumber \\ 
& = 1 + \frac{1}{3 i} \bigg( e^{i \pi/2} H_{0,1}^{1,0}\left[- e^{-i \pi/6} \sqrt{3} x \bigg| \begin{array}{c}-\\ (0,1)
\end{array}\right]  \nonumber \\
&- e^{-i \pi/2} H_{0,1}^{1,0}\left[- e^{i \pi/6} \sqrt{3} x \bigg| \begin{array}{c}
-\\
(0,1)
\end{array}\right] \bigg)
\end{align}
Using (\ref{exph}) on (\ref{hypnbai6-1}) leads to:
\begin{equation} \label{hypnbai6-2}
\frac{1}{3}+ \frac{2 \pi}{3}  
H_{1,2}^{1,0}\left[-\sqrt{3} x \bigg| \begin{array}{c}
(\frac{1}{2},\frac{1}{6})\\
(0,1),(\frac{1}{2},\frac{1}{6})
\end{array}\right] = \frac{1}{3} + \frac{2}{3} e^{3x/2} cos \left(\frac{\sqrt{3} x}{2} \right)
\end{equation}

\noindent which is an alternate form of (\ref{hypnbai6}) in terms of trigonometrical functions.

\subsubsection{Case when $\delta = 1$}

When $\delta = 1$, (\ref{hypnbai7}) and (\ref{eq:HfuncIDAKR}) provide:

\begin{align} \label{hypnbai7-1}
\frac{2 \pi}{\sqrt{3}} H_{2,3}^{1,1}& \left[-\sqrt{3} x \bigg| \begin{array}{c} (0,1),(\frac{2}{3},\frac{1}{6}) \\ (0,1),(-1,1),(\frac{2}{3},\frac{1}{6})\end{array}\right] \nonumber \\ 
&= \frac{1}{\sqrt{3}i} \bigg( e^{2i \pi/3} H_{1,2}^{1,1}\left[-e^{-i \pi/6} \sqrt{3} x \bigg| \begin{array}{c}(0,1)\\(0,1),(-1,1)\end{array}\right]  \nonumber\\
& - e^{-2i \pi/3} H_{1,2}^{1,1}\left[-e^{i \pi/6} \sqrt{3} x \bigg| \begin{array}{c}
(0,1)\\ (0,1),(-1,1)\end{array}\right] \bigg)
\end{align}

By combining (\ref{h2}) and (\ref{hypnbai7-1}), the following result is obtained:

\begin{equation} \label{hypnbai7-2}
\frac{2 \pi}{\sqrt{3}}  
H_{2,3}^{1,1}\left[-\sqrt{3} x \bigg| \begin{array}{c}
(0,1),(\frac{2}{3},\frac{1}{6})\\
(0,1),(-1,1),(\frac{2}{3},\frac{1}{6})
\end{array}\right] = \frac{2}{3x} \left(\frac{1}{2} +e^{\frac{3x}{2}} sin \left( \frac{5 \pi}{6} - \frac{\sqrt{3} x}{2}\right) \right) 
\end{equation}

\noindent which is an alternate form for (\ref{hypnbai7}).

\subsubsection{Case when $\delta = 2$}

In this case, the combination of (\ref{hypnbai8}) and (\ref{eq:HfuncIDAKR}) leads to:

\begin{align} \label{hypnbai8-1}
4 \pi  H_{3,4}^{1,2}& \left[-\sqrt{3} x \bigg| \begin{array}{c} (0,1),(-1,1),(\frac{5}{6},\frac{1}{6})\\ (0,1),(-1,1),(-2,1),(\frac{5}{6},\frac{1}{6}) \end{array}\right] \nonumber \\ 
 & = \frac{2}{i} \bigg( e^{5 i \pi/6} H_{2,3}^{1,2}\left[-e^{-i \pi/6}\sqrt{3} x \bigg| \begin{array}{c}(0,1),(-1,1)\\ (0,1),(-1,1),(-2,1)
\end{array}\right]  \nonumber \\ 
& - e^{-5 i \pi/6} H_{2,3}^{1,2}\left[-e^{i \pi/6}\sqrt{3} x \bigg| \begin{array}{c}
(0,1),(-1,1)\\ (0,1),(-1,1),(-2,1) \end{array}\right] \bigg)
\end{align}

It is clear from (\ref{h3}) that (\ref{hypnbai8-1}) may be rewritten as:

\begin{equation} \label{hypnbai8-2}
4 \pi  
H_{3,4}^{1,2}\left[-\sqrt{3} x \bigg| \begin{array}{c}
(0,1),(-1,1),(\frac{5}{6},\frac{1}{6})\\
(0,1),(-1,1),(-2,1),(\frac{5}{6},\frac{1}{6})
\end{array}\right] = \frac{2}{3 x^2} \left(1-2e^{\frac{3x}{2}} cos \left( \frac{\pi}{3} + \frac{\sqrt{3}x}{2} \right) \right)
\end{equation}

\noindent giving and alternate form for (\ref{hypnbai8}) in terms of elementary functions.

\subsection{When $\gamma=2$ and $\beta =1$ in (\ref{refhyp})}
The cases where $\beta =1$, $\gamma=3$ ($\delta = 0, 1$ and $2$) and $\gamma=4$ ($\delta = 0, 1, 2$ and $3$) have been explored in \cite{HMS2012}. In the latter, by taking $\Omega_i = 1$,$\forall i \in \mathbb{N}$, the authors provide some general results for summations similar to the ones obtained the earlier subsections.

In this subsection, we obtain two results not given in \cite{HMS2012}, by taking  $\gamma=2$ and $\beta =1$ in (\ref{refhyp}) and (\ref{refhyp2}) utilizing (\ref{eq:HfuncID}).

For example, the result (2.1) for $\Omega_i = 1$,$\forall i \in \mathbb{N}$ taken in \cite{HMS2012} is written in the following form:

\begin{eqnarray} \label{nr1}
e^{x}  
H_{1,2}^{1,1}\left[- x^3 \bigg| \begin{array}{c}
(0,1)\\
(0,1),(0,3)
\end{array}\right]& = & e^{x} {}_{0} F_{2} \left[ \begin{array}{c}
- \\
\frac{1}{3}, \frac{2}{3}
\end{array} \; ; \frac{x^{3}}{27} \right] \\ \nonumber
& = & \frac{2}{3} \sum_{m=0}^{\infty} \frac{x^m}{m!} \left[ 2^{m-1} + cos \left(\frac{m \pi}{3} \right) \right] \\ \nonumber
& = & \frac{e^{2x}}{3} + \frac{2 \pi}{3} H_{1,2}^{1,0}\left[- x \bigg| \begin{array}{c}
(\frac{1}{2},\frac{1}{3})\\
(0,1),(\frac{1}{2},\frac{1}{3})
\end{array}\right] \\ \nonumber
& = & \frac{e^{2x}}{3} + \frac{2}{3}e^{\frac{x}{2}} cos \left(\frac{\sqrt{3} x}{2} \right)
\end{eqnarray}

The first equality in (\ref{nr1}) is easily obtained from  (\ref{refhyp}), (\ref{refhyp2}) and (2.1) for $\Omega_i = 1$,$\forall i \in \mathbb{N}$ taken in \cite{HMS2012}. The third and fourth inequalities are obtained below:

We start with:

\begin{eqnarray} \label{nr2}
\frac{2}{3} \sum_{m=0}^{\infty} \frac{x^m}{m!} \left[ 2^{m-1} + cos \left(\frac{m \pi}{3} \right) \right] 
& = & \frac{1}{3} \sum_{m=0}^{\infty} \frac{(2x)^m}{m!} + \frac{2}{3} \sum_{m=0}^{\infty} \frac{x^m}{m!} cos \left(\frac{m \pi}{3} \right) \\ \nonumber
& = & \frac{e^{2x}}{3} + \frac{2 \pi}{3} \sum_{m=0}^{\infty} \frac{x^m}{m! \Gamma \left(\frac{1}{2} + \frac{m}{3} \right) \Gamma \left(\frac{1}{2} - \frac{m}{3} \right)} \\ \nonumber
& = & \frac{e^{2x}}{3} + \frac{2 \pi}{3} H_{1,2}^{1,0}\left[- x \bigg| \begin{array}{c}
(\frac{1}{2},\frac{1}{3})\\
(0,1),(\frac{1}{2},\frac{1}{3})
\end{array}\right]
\end{eqnarray}

\noindent using the series representation in (\ref{hcomput1}). On the other hand,
\begin{equation}
\frac{e^{2x}}{3} + \frac{2 \pi}{3} H_{1,2}^{1,0}\left[- x \bigg| \begin{array}{c}
(\frac{1}{2},\frac{1}{3})\\
(0,1),(\frac{1}{2},\frac{1}{3})
\end{array}\right] = \frac{e^{2x}}{3} + \frac{2 \pi}{3} \frac{e^{x cos (\pi/3)}}{\pi} sin \left(\frac{\pi}{2} -x sin \left(\frac{\pi}{3} \right)\right)
\end{equation}

\noindent by using \cite{AKR1981} as done in \S 5.1.1. Finally:

\begin{equation}
\frac{e^{2x}}{3} + \frac{2 \pi}{3} H_{1,2}^{1,0}\left[- x \bigg| \begin{array}{c}
(\frac{1}{2},\frac{1}{3})\\
(0,1),(\frac{1}{2},\frac{1}{3})
\end{array}\right] = \frac{e^{2x}}{3} + \frac{2}{3}e^{\frac{x}{2}} cos \left(\frac{\sqrt{3} x}{2} \right)
\end{equation}

The other six results given in \cite{HMS2012} with $\Omega_i = 1$,$\forall i \in \mathbb{N}$ can be written in terms of H-functions by using (\ref{hcomput1}) or (\ref{hcomput2}) and in elementary functions utilizing \cite{AKR1981, Brychkov}. These results are too complicated to be included here.

\subsubsection{Case when $\delta = 0$}

The case where $\delta = 0$ in (\ref{refhyp2}) implies:

\begin{equation} \label{refhyp4}
\sum_{n=0}^{\infty} \frac{x^{\alpha n}}{(2 n)!} =  
H_{1,2}^{1,1}\left[- x^{\alpha}|\begin{array}{c}
(0,1)\\
(0,1),(0,2)
\end{array}\right] 
\end{equation}

On the other hand, the result in (\ref{eq:HfuncID}) provides:

\begin{eqnarray}
H_{1,2}^{1,1}\left[- x^{\alpha}|\begin{array}{c}
(0,1)\\
(0,1),(0,2)
\end{array}\right]  = H_{1,2}^{1,1}\left[- e^{- i \pi} x^{\alpha}\;\bigg|\begin{array}{c}
(0,2)\\
(0,2), (0,2)
\end{array}\right] \nonumber \\
+  H_{1,2}^{1,1}\left[- e^{i \pi} x^{\alpha}\;\bigg|\begin{array}{c}
(0,2)\\
(0,2), (0,2)
\end{array}\right],\label{eq:HfuncID2}
\end{eqnarray}

Also, by considering the identity (\ref{hexp}) with $k=1/2$, (\ref{eq:HfuncID2})becomes:

\begin{align}
H_{1,2}^{1,1}\left[- x^{\alpha}|\begin{array}{c} (0,1)\\ (0,1),(0,2)
\end{array}\right]  & = \frac{1}{2} H_{1,2}^{1,1}\left[x^{\alpha/2}\;\bigg|\begin{array}{c}
(0,1)\\ (0,1), (0,1) \end{array}\right] \nonumber \\
& \;  + \frac{1}{2} H_{1,2}^{1,1}\left[- x^{\alpha/2}\;\bigg|\begin{array}{c}
(0,1)\\ (0,1), (0,1) \end{array}\right],\label{eq:HfuncID3}
\end{align}

By further manipulating the contour integral representation of the H-functions in \ref{eq:HfuncID3}, it is clear from (\ref{exph}) that (\ref{eq:HfuncID3}) turns to \cite{Brychkov}:

\begin{equation}
H_{1,2}^{1,1}\left[- x^{\alpha}|\begin{array}{c}
(0,1)\\
(0,1),(0,2)
\end{array}\right]   = \frac{1}{2} e^{x^{\alpha/2}} 
+ \frac{1}{2} e^{-x^{\alpha/2}} 
 = cosh(x^{\alpha/2}),\label{eq:HfuncID4}
\end{equation}

\subsubsection{Case when $\delta = 1$}

When $\delta = 1$:

\begin{equation} \label{refhyp4-2}
\sum_{n=0}^{\infty} \frac{x^{\alpha n}}{(2 n + 1)!} =
H_{1,2}^{1,1}\left[- x^{\alpha}|\begin{array}{c}
(0,1)\\
(0,1),(-1,2)
\end{array}\right] 
\end{equation}

On the other hand, the result in (\ref{eq:HfuncID}) provides:

\begin{eqnarray}
H_{1,2}^{1,1}\left[- x^{\alpha}|\begin{array}{c}
(0,1)\\
(0,1),(-1,2)
\end{array}\right]  = H_{1,2}^{1,1}\left[- e^{- i \pi} x^{\alpha}\;\bigg|\begin{array}{c}
(0,2)\\
(0,2), (-1,2)
\end{array}\right] \nonumber \\
+  H_{1,2}^{1,1}\left[- e^{i \pi} x^{\alpha}\;\bigg|\begin{array}{c}
(0,2)\\
(0,2), (-1,2)
\end{array}\right],\label{eq:HfuncID2-2}
\end{eqnarray}

Also, by considering the identity (\ref{hexp}) with $k=1/2$, (\ref{eq:HfuncID2-2})becomes:

\begin{align}
H_{1,2}^{1,1}\left[- x^{\alpha}|\begin{array}{c}
(0,1)\\
(0,1),(-1,2)
\end{array}\right] & = \frac{1}{2} H_{1,2}^{1,1}\left[x^{\alpha/2}\;\bigg|\begin{array}{c}
(0,1)\\
(0,1), (-1,1)
\end{array}\right] \nonumber \\
& \; + \frac{1}{2} H_{1,2}^{1,1}\left[- x^{\alpha/2}\;\bigg|\begin{array}{c}
(0,1)\\
(0,1), (-1,1)
\end{array}\right],\label{eq:HfuncID3-2}
\end{align}

By using the result (\ref{h2}), one may see that (\ref{eq:HfuncID3-2}) turns to \cite{Brychkov}:

\begin{align}
H_{1,2}^{1,1}\left[- x^{\alpha}|\begin{array}{c}
(0,1)\\
(0,1),(-1,2)
\end{array}\right] & = \frac{(1-e^{-x^{\alpha/2}}) x^{-\alpha/2}}{2} 
+ \frac{(e^{x^{\alpha/2}}-1) x^{-\alpha/2}}{2} \nonumber \\
& = x^{-\alpha/2} sinh(x^{\alpha/2}),\label{eq:HfuncID4-2}
\end{align}

\section{Conclusions}

A new identity which enables one to split certain H-functions into the sum of two othe H-functions has been derived. This new formula has been applied to simplify the summation of hypergeometric series. By using this identity, new relations between H-functions have been established.

\section*{Conflict of interests}
No conflict of interests was reported by the authors.

\section*{Acknowledgments}
P. N. Rathie thanks the Coordination for the Improvement of Higher Level Personnel (CAPES) for supporting his Senior National Visiting Professorship. Also, L. C. de S. M. Ozelim thanks the Brazilian
Research Council (CNPq) for funding his postdoctoral fellowship at University of Brasilia.

\end{document}